\documentclass[11pt]{amsart}

\usepackage[a4paper,margin=33mm]{geometry}
\usepackage[T1]{fontenc}
\usepackage{newtxtext,newtxmath}
\usepackage{microtype}
\usepackage{mathtools}
\usepackage{enumitem}
\usepackage{aliascnt}
\usepackage[dvipsnames]{xcolor}
\usepackage[colorlinks=true,linkcolor=blue!45!black,citecolor=green!35!black,
  urlcolor=blue!45!black]{hyperref}
\usepackage[nameinlink,noabbrev]{cleveref}

\numberwithin{equation}{section}
\setlist[enumerate,1]{label=\textup{(\roman*)},ref=\textup{\roman*},leftmargin=2.2em}

\newtheorem{theorem}{Theorem}[section]
\crefname{theorem}{Theorem}{Theorems}
\Crefname{theorem}{Theorem}{Theorems}

\newaliascnt{problem}{theorem}

\aliascntresetthe{problem}
\crefname{problem}{problem}{problems}
\Crefname{problem}{Problem}{Problems}

\newaliascnt{lemma}{theorem}
\newtheorem{lemma}[lemma]{Lemma}
\aliascntresetthe{lemma}
\crefname{lemma}{Lemma}{lemmas}
\Crefname{lemma}{Lemma}{Lemmas}

\newaliascnt{corollary}{theorem}
\newtheorem{corollary}[corollary]{Corollary}
\aliascntresetthe{corollary}
\crefname{corollary}{Corollary}{Corollaries}
\Crefname{corollary}{Corollary}{Corollaries}

\newaliascnt{definition}{theorem}

\aliascntresetthe{definition}
\crefname{definition}{Definition}{definitions}
\Crefname{definition}{Definition}{Definitions}

\newaliascnt{proposition}{theorem}
\newtheorem{proposition}[proposition]{Proposition}
\aliascntresetthe{proposition}
\crefname{proposition}{Proposition}{Propositions}
\Crefname{proposition}{Proposition}{Propositions}

\newaliascnt{result}{theorem}

\aliascntresetthe{result}
\crefname{result}{result}{results}
\Crefname{result}{Result}{Results}

\newaliascnt{example}{theorem}
\newtheorem{example}[example]{Example}
\aliascntresetthe{example}
\crefname{example}{example}{examples}
\Crefname{example}{Example}{Examples}

\newaliascnt{conjecture}{theorem}

\aliascntresetthe{conjecture}
\crefname{conjecture}{conjecture}{conjectures}
\Crefname{conjecture}{Conjecture}{Conjectures}

\newaliascnt{observ}{theorem}

\aliascntresetthe{observ}
\crefname{observ}{observation}{observations}
\Crefname{observ}{Observation}{Observations}

\newaliascnt{remark}{theorem}
\newtheorem{remark}[remark]{Remark}
\aliascntresetthe{remark}
\crefname{remark}{remark}{remarks}
\Crefname{remark}{Remark}{Remarks}

\newcommand{\F}{\mathbb{F}}
\newcommand{\C}{\mathbb{C}}
\newcommand{\Qzero}{Q_0}
\newcommand{\Qone}{Q_1}
\newcommand{\supp}{\operatorname{supp}}
\newcommand{\restr}{\mathbin{|}}

\makeatletter
\@namedef{subjclassname@2020}{%
  \textup{2020} Mathematics Subject Classification}
\makeatother

\title[Fourier bases for square-order Paley graphs]
  {Fourier Bases from Optimal Eigenfunctions\\
   of Paley Graphs of Square Order}

\author{Bart De Bruyn$^{2}$}
\author{Sergey Goryainov$^1$}
\author{Yilin Xie$^1$}
\author{Weihao Yan$^{1,2}$}
\author{Ferdinando Zullo$^{1,3}$}
\address{$^1$School of Mathematical Sciences, Hebei International Joint Research Center for Mathematics and Interdisciplinary Science, Hebei Key Laboratory of Computational Mathematics and Applications, Hebei Workstation for Foreign Academicians, Hebei Normal University, Shijiazhuang 050024,
P. R. China}
\address{$^2$Department of Mathematics, Computer Science and Statistics, Ghent University, Ghent, Belgium}
\address{$^3$Dipartimento di Matematica e Fisica, Universit\`a degli Studi della Campania ``Luigi Vanvitelli'', I--\,81100 Caserta, Italy}

\email{Bart.DeBruyn@ugent.be, sergey.goryainov3@gmail.com, ferdinando.zullo@unicampania.it}

\date{}

\begin{document}

\begin{abstract}
Let $q$ be an odd prime power and consider the Paley graph $P(q^2)$ on the finite field $\F_{q^2}$. A well-known optimal eigenfunction of this graph is supported on the norm-one subgroup, where it takes the values $1$ and $-1$ on the two cosets of its subgroup of index two. Its support has the minimum possible size $q+1$.
Starting from this eigenfunction, we construct a family of multiplicative translates and prove that they form a basis of the space of complex-valued functions on $\F_{q^2}^{\times}$. After extending these functions to the whole field, the translates indexed by squares form a basis of one nonprincipal eigenspace of the Paley graph, while those indexed by nonsquares form a basis of the other. In particular, both nonprincipal eigenspaces admit bases consisting entirely of optimal eigenfunctions.
The main ingredient is a multiplicative Fourier analysis of the shifted eigenfunction. We determine the moduli of all its Fourier coefficients and, in particular, show that none of them vanishes.
\end{abstract}

\maketitle

\pagestyle{plain}

\section{Introduction}

Let $q$ be an odd prime power.  The Paley graph $P(q^2)$ is the Cayley
graph on the additive group of $F=\F_{q^2}$ in which distinct vertices
$x$ and $y$ are adjacent if $x-y$ is a square in $F^\times$.  Paley
graphs arose from Paley's construction of conference matrices
\cite{Paley1933}. They have wide-ranging applications across combinatorics, 
including in Ramsey theory, design theory, and coding theory. Many of their interesting properties 
stem from their pseudo-random structure: for instance, every finite graph occurs 
as an induced subgraph of some sufficiently large Paley graph (see e.g. 
\cite[Theorem 6.19]{Nica}).

Paley graphs also play a prominent role in spectral graph theory and algebraic 
combinatorics in general. Standard accounts of their spectral and strongly
regular properties may be found in \cite{BrouwerVanMaldeghem2022,Jones2020}.
In the square-order case, $P(q^2)$ has valency $(q^2-1)/2$, nonprincipal
eigenvalues
\[
  r=\frac{q-1}{2},\qquad s=-\frac{q+1}{2},
\]
and both nonprincipal eigenspaces have dimension $(q^2-1)/2$.
We regard these eigenspaces as complex vector spaces unless stated
otherwise.  Since all the basis functions constructed below are
real-valued, the same families will also be bases of the corresponding
real eigenspaces.

The square-order graphs sit naturally in finite geometry.  Viewing
$F$ as a two-dimensional vector space over $\F_q$ identifies it with
$\operatorname{AG}(2,q)$.  The affine lines then split into square and
nonsquare types: the former are cliques of size $q$ and the latter are
cocliques of size $q$.  These sets attain the Delsarte--Hoffman bound
\cite{Delsarte1973}.
Motivated by a problem about generalized quadratic-residue codes
\cite{VanLintMacWilliams1978}, Blokhuis proved that they are the only
maximum cliques and cocliques \cite{Blokhuis1984}.  Sziklai extended the
underlying power-difference result beyond squares \cite{Sziklai1999};
recent developments and further extensions are discussed in
\cite{AsgarliYip2022,BrouwerEtAl2026,Yip2026}.

Maximal cliques below the Delsarte bound display a different geometry.
Baker, Ebert, Hemmeter, and Woldar constructed maximal cliques of size
$(q+1)/2$ for $q\equiv1\pmod4$ and $(q+3)/2$ for
$q\equiv3\pmod4$ \cite{BakerEtAl1996}.  Related constructions appear in
\cite{KiermaierKurz2009}.  Goryainov, Kabanov, Shalaginov, and
Valyuzhenich connected two such configurations to eigenfunctions with
support $q+1$ \cite{GoryainovEtAl2018}; this is the minimum allowed by
the weight-distribution bound.  The minimum-support problem for graph
eigenfunctions is surveyed in \cite{SotnikovaValyuzhenich2021}.
Goryainov, Masley, and Shalaginov later exhibited a linear-fractional
correspondence between the two clique constructions
\cite{GoryainovMasleyShalaginov2022}.  More recently, Brouwer,
Goryainov, Shalaginov, and Yip placed these constructions in the wider
setting of Desarguesian nets and Peisert graphs
\cite{BrouwerEtAl2026}.
The related eigenfunction construction has also been extended to
generalized Paley graphs of square order
\cite{GoryainovShalaginovYip2023}.

The present paper turns from the maximality of the underlying cliques
to the span of their associated optimal eigenfunctions.  Let
\[
  Q=\{u\in F^\times:N_{F/\F_q}(u)=1\}
  =\ker N_{F/\F_q}
\]
be the norm-one subgroup of $F^\times$.  It is cyclic of order $q+1$
and has a unique subgroup $\Qzero$ of index $2$.  Put
$\Qone=Q\setminus\Qzero$ and define $f\colon F\to\C$ by
\[
  f(x)=
  \begin{cases}
    1,&x\in\Qzero,\\
    -1,&x\in\Qone,\\
    0,&x\notin Q.
  \end{cases}
\]
We call a nonzero $\lambda$-eigenfunction \emph{optimal} if its support
has minimum cardinality among all nonzero $\lambda$-eigenfunctions.
The work \cite{GoryainovEtAl2018} cited above constructs such eigenfunctions and relates
the underlying clique configurations.  Our contribution concerns a
different question: their linear span. We show first that all multiplicative translates $h_z(t) = f(1 + \frac{t}{z})$, $z \in F^\times$, of the shifted function $h(t)=f(1+t)$ form a
basis of $\C^{F^\times}$.  After extending these functions back to
$F$, the square and nonsquare multiplicative parameters $z$ give bases of the
two nonprincipal eigenspaces.

\begin{theorem}[Main theorem]
\label{thm:main}
Let $q$ be an odd prime power, let $F=\F_{q^2}$, and let $P(q^2)$ be
the Paley graph on $F$.  Define $Q=\ker N_{F/\F_q}$, let $\Qzero$ be
the subgroup of $Q$ of order $(q+1)/2$, and put
$\Qone=Q\setminus\Qzero$.  Let $f\colon F\to\C$ be given by
\[
  f(x)=
  \begin{cases}
    1,&x\in\Qzero,\\
    -1,&x\in\Qone,\\
    0,&x\notin Q.
  \end{cases}
\]
For $z\in F^\times$, define
\[
  h_z(x)=f\left(1+\frac{x}{z}\right),
  \qquad x\in F^\times.
\]
Then $\{ h_z:z\in F^\times \}$ is a basis of $\C^{F^\times}$.

Let $D=(F^\times)^2$ and extend $h_z$ to $F$ by
\[
  \widetilde h_z(x)=f\left(1+\frac{x}{z}\right),
  \qquad x\in F.
\]
The family $\{\widetilde h_z:z\in D\}$ is a basis of the same
nonprincipal eigenspace as $f$, while
$\{\widetilde h_z:z\in F^\times\setminus D\}$ is a basis of the other
nonprincipal eigenspace.  
Every $\widetilde h_z$ is an optimal eigenfunction with support of
size $q+1$.  The full family
$\{\widetilde h_z:z\in F^\times\}$ is a basis of the direct sum of the
two nonprincipal eigenspaces, namely of $\mathbf1^\perp$.  Since these
functions are real-valued, all the same assertions hold over
$\mathbb R$.
\end{theorem}

\medskip As a consequence of Theorem \ref{thm:main} we will obtain characterizations for the eigenfunctions of the nonprincipal eigenvalues. These characterizations invoke the following properties for a function $g: F \to \mathbb{R}$:
\begin{enumerate}
\item[(P1)] $\sum_{x \in F} g(x) = 0$;
\item[(P2)] $\sum_{x \in z(\Qzero-1)} g(x) = \sum_{x \in z(\Qone-1)} g(x)$ for every square $z \in F^\times$;
\item[(P3)] $\sum_{x \in z(\Qzero-1)} g(x) = \sum_{x \in z(\Qone-1)} g(x)$ for every nonsquare $z \in F^\times$.
\end{enumerate}
Here, $z(A-1) := \{ z(\alpha-1) \, : \,  \alpha \in A \}$ for every $A \in \{ Q,\Qzero,\Qone \}$ and every $z \in F^\times$.

\begin{corollary} \label{co:main}
\begin{itemize}
\item If $\lambda \equiv 1 \pmod{4}$, then a nonzero function $g: F^\times \to \mathbb{R}$ is an eigenfunction for the eigenvalue $\frac{q-1}{2}$ of $P(q^2)$ if and only if (P1) and (P2) hold, and it is an eigenfunction for the eigenvalue $\frac{-q-1}{2}$ if and only if (P1) and (P3) hold. 
\item If $\lambda \equiv 3 \pmod{4}$, then a nonzero function $g: F^\times \to \mathbb{R}$ is an eigenfunction for the eigenvalue $\frac{q-1}{2}$ of $P(q^2)$ if and only if (P1) and (P3) hold, and it is an eigenfunction for the eigenvalue $\frac{-q-1}{2}$ if and only if (P1) and (P2) hold. 
\end{itemize}
\end{corollary}

\begin{remark}
For a function $g: F \to \mathbb{R}$ to belong to the eigenspace corresponding to $\lambda \in \{ \frac{q-1}{2},\frac{-q-1}{2} \}$, a certain collection of $\frac{q^2+1}{2}$ conditions must thus be satisfied. These $\frac{q^2+1}{2}$ conditions are really necessary, as omitting only one of them already results in counter examples.
\end{remark}

\medskip The key step in the proof of our results is independent of the additive Cayley-graph structure. We calculate the modulus of every multiplicative Fourier coefficient of $h$ and show that none vanishes.  It follows that all $q^2-1$ multiplicative translates of $h$ are linearly independent. The eigenspace statements then follow by separating the square and nonsquare parameters $z$.

In Section \ref{sec:geometry} we collect the finite-field and spectral facts
needed below.  The Fourier calculation is carried out in
Section \ref{sec:fourier}; immediately after the Fourier-basis theorem we
derive the two eigenspace bases by separating the square and nonsquare
values of $z$.  In Section \ref{sec:examples} we work out these two choices
explicitly for $q=5$ and $q=3$.

\section{The norm-one oval and its eigenfunction}
\label{sec:geometry}

\subsection{The norm-one group}

Throughout the paper, $F=\F_{q^2}$ and
\[
  N=N_{F/\F_q}\colon F^\times\longrightarrow\F_q^\times,
  \qquad N(x)=x^{q+1}.
\]
The group
\[
  Q=\ker N
\]
is cyclic of order $q+1$.  Let $\Qzero$ be its unique subgroup of
order $(q+1)/2$, put $\Qone=Q\setminus\Qzero$, and denote by
\[
  \chi\colon Q\longrightarrow\{\pm1\}
\]
the quadratic character of $Q$, whose kernel is $\Qzero$.  We extend
$\chi$ by zero to a function $f\colon F\to\C$:
\begin{equation}
  f(x)=
  \begin{cases}
    \chi(x),&x\in Q,\\
    0,&x\notin Q.
  \end{cases}
  \label{eq:def-f}
\end{equation}

\subsubsection*{The different roles of $Q$}

Although $Q$ appears below in algebraic, geometric, and Fourier
descriptions, the symbol always denotes the same set: the norm-one
subgroup of $F^\times$.  We now make the relations between these
descriptions explicit.

First, $Q$ is the norm-one subgroup of $F^\times$:
\[
  Q=\{x\in F^\times:N(x)=1\}.
\]
Equivalently, it is the kernel of the surjective norm map
\[
  N\colon F^\times\longrightarrow\F_q^\times.
\]
In particular,
\[
  F^\times/Q\cong\F_q^\times.
\]
This quotient is not the one used in the Fourier calculation, but it
records the primary definition of $Q$ as a kernel.

Second, the map
\[
  F^\times\longrightarrow Q,
  \qquad r\longmapsto r^{q-1},
\]
is surjective and has kernel $\F_q^\times$.  It therefore induces an
isomorphism
\begin{equation}
  F^\times/\F_q^\times\xrightarrow{\ \sim\ }Q,
  \qquad r\F_q^\times\longmapsto r^{q-1}.
  \label{eq:quotient-to-Q}
\end{equation}
Thus $Q$ is not redefined as a quotient group: rather, the norm-one
subgroup is canonically identified with this quotient.  This is the
identification used in Section \ref{sec:fourier} to transport characters
that are trivial on $\F_q^\times$ to characters of $Q$.

Third, when $F$ is regarded as a two-dimensional vector space over
$\F_q$, its elements are points of $\operatorname{AG}(2,q)$.  Under
this identification, the same subset $Q\subset F$ is a nonsingular
affine conic; the coordinate equation is given below.  The conic is
the support of the eigenfunction $f$ and underlies the geometric
secant-and-tangent proof of \cref{prop:optimal} given in
\cite{GoryainovEtAl2018}.

Finally, $\Qzero$ is a subgroup of index $2$ in $Q$, and
\begin{equation}
  Q/\Qzero\cong\{\pm1\}.
  \label{eq:Q-mod-Q0}
\end{equation}
The character $\chi$ is the quotient character associated with
\eqref{eq:Q-mod-Q0}: it equals $1$ on $\Qzero$ and $-1$ on the other
coset $\Qone$.  Accordingly, the roles used in the paper may be
summarized as follows:
\[
\begin{array}{c@{\qquad}l}
Q\le F^\times
  &\text{the cyclic norm-one subgroup},\\
Q\subset\operatorname{AG}(2,q)
  &\text{the norm-one affine conic},\\
Q\cong F^\times/\F_q^\times
  &\text{the quotient model used to transport characters},\\
Q/\Qzero\cong\{\pm1\}
  &\text{the quotient that defines $\chi$ and hence $f$}.
\end{array}
\]

The following elementary facts will be used repeatedly.

\begin{lemma}
\label{lem:basic-Q}
The norm-one group satisfies
\begin{enumerate}
  \item $u^q=u^{-1}$ for every $u\in Q$;
  \item $\F_q^\times\cap Q=\{\pm1\}$;
  \item every element of $Q$ and every element of $\F_q^\times$ is a
  square in $F^\times$;
  \item
  \[
    \chi(-1)=
    \begin{cases}
      -1,&q\equiv1\pmod4,\\
      1,&q\equiv3\pmod4.
    \end{cases}
  \]
\end{enumerate}
\end{lemma}

\begin{proof}
If $u\in Q$, then $u^{q+1}=1$, and hence $u^q=u^{-1}$.  Now let
$x\in\F_q^\times\cap Q$.  Since $x^q=x$ and $x^{q+1}=1$, we have
$x^2=1$.  Thus $x=\pm1$, and the reverse inclusion is immediate.

Let $F^\times=\langle g\rangle$.  Then
\[
  Q=\langle g^{q-1}\rangle,
  \qquad
  \F_q^\times=\langle g^{q+1}\rangle.
\]
Both $q-1$ and $q+1$ are even, so both subgroups lie in
$(F^\times)^2$.

The element $-1$ has order $2$ in the cyclic group $Q$.  It belongs to
$\Qzero$ precisely when $2$ divides $|\Qzero|=(q+1)/2$.  This is
equivalent to $q\equiv3\pmod4$ and gives the displayed values of
$\chi(-1)$.
\end{proof}

\medskip
\begin{lemma}
\label{lem:Q-equivariance}
For all $b\in Q$ and $y\in F$,
\[
  f(by)=\chi(b)f(y).
\]
\end{lemma}

\begin{proof}
If $y\notin Q$, then $by\notin Q$, and both sides vanish.  If
$y\in Q$, then
\[
  f(by)=\chi(by)=\chi(b)\chi(y)=\chi(b)f(y).
\]
\end{proof}

For a geometric interpretation, choose a nonsquare $d$ in $\F_q^\times$
and write $F=\F_q(\alpha)$ with $\alpha^2=d$.  Every element of $F$ has
a unique form $X+Y\alpha$ with $X,Y \in \mathbb{F}_q$, and
\[
  N(X+Y\alpha)=X^2-dY^2.
\]
Thus $Q$ is the affine conic $X^2-dY^2=1$ in
$\operatorname{AG}(2,q)$.

\subsection{Optimal eigenfunctions}

Paley graphs are examples of primitive strongly regular graphs. For such graphs, the weight-distribution bound takes a simple form and a combinatorial interpretation in case of tightness can be given, as indicated in the following two lemmas. These results were derived in \cite[Corollary 2 and Lemma 3]{Go-Pa} as special cases of a more general result that can be found in \cite{Kr-Mo-Po}.

\begin{lemma} \label{newa}
Let $G$ be a primitive strongly regular graph with nonprincipal eigenvalues $s<0<r$. Then an eigenfunction of $G$ corresponding to the eigenvalue $r$ has at least $2(r+1)$ nonzeroes, and an eigenfunction corresponding to the eigenvalue $s$ has at least $-2s$ nonzeroes.
\end{lemma}

\begin{lemma} \label{newb}
Let $G$ be a primitive strongly regular graph with nonprincipal eigenvalues $s,r$, where $s < 0 < r$. Then the following statements hold.
\begin{enumerate}
\item[$(1)$] For an $s$-eigenfunction $f$, if the cardinality of support of $f$ meets the weight-distribution bound, then there exists an induced complete bipartite subgraph in $G$ with parts $T_0$ and $T_1$ of size $-s$. Moreover, up to multiplication by a constant, $f$ has value $1$ on the vertices of $T_0$ and value $-1$ on the vertices of $T_1$.
\item[$(2)$] For an $r$-eigenfunction $f$, if the cardinality of support of $f$ meets the weight-distribution bound, then there exists an induced pair of isolated cliques $T_0$ and $T_1$ in $G$ of size $-\overline{s} = -(-1-r) = 1+r$. Moreover, up to multiplication by a constant, $f$ has value 1 on the vertices of $T_0$ and value $-1$ on the vertices of $T_1$.
\end{enumerate}
\end{lemma}

\medskip \noindent For the Paley graph $P(q^2)$, the weight-distribution bound is indeed tight. The following proposition, taken from Goryainov, Kabanov, Shalaginov and Valyuzhenich \cite[Lemma~8 and Theorem~2]{GoryainovEtAl2018}, makes the claims of Lemmas \ref{newa} and \ref{newb} more explicit.

\begin{proposition}
\label{prop:optimal}
The function $f$ defined in \eqref{eq:def-f} is an eigenfunction of
$P(q^2)$ with support $Q$, and hence
\[
  |\supp(f)|=q+1.
\]
More precisely:
\begin{enumerate}
  \item if $q\equiv1\pmod4$, then $f$ has eigenvalue
  $-(q+1)/2$, and the subgraph induced on $Q$ is the complete
  bipartite graph with parts $\Qzero$ and $\Qone$;
  \item if $q\equiv3\pmod4$, then $f$ has eigenvalue
  $(q-1)/2$, and the subgraph induced on $Q$ is the disjoint union of
  the cliques on $\Qzero$ and $\Qone$.
\end{enumerate}
In both cases, $f$ has minimum support among the eigenfunctions
belonging to the corresponding eigenvalue.
\end{proposition}

This result, which will be useful later, has a simple geometric interpretation.  The function $f$
takes the value $1$ on $\Qzero$, the value $-1$ on $\Qone$, and
vanishes outside the conic $Q$.  When $q\equiv1\pmod4$, adjacency
inside $Q$ exchanges $\Qzero$ and $\Qone$, which accounts for the
negative eigenvalue.  When $q\equiv3\pmod4$, adjacency preserves the
two parts, which accounts for the positive eigenvalue.  The cited
result proves that the eigenfunction equation also holds at the
vertices outside $Q$.  Moreover, the support size $q+1$ attains the
weight-distribution bound; see also
\cite{SotnikovaValyuzhenich2021}.

\section{Nonvanishing Fourier coefficients} \label{sec:fourier}

Define
\[
  h\colon F^\times\longrightarrow\C,
  \qquad h(t)=f(1+t).
\]
For $z\in F^\times$, let $h_z(x)=h(x/z)$.  The relevance of these
multiplicative translates is given by the following observation.

\begin{lemma}
\label{lem:reduction}
For $a\in\F_q^\times$ and $b\in Q$, put $z=b/a$.  Then
\[
  f(ax+b)=\chi(b)h_z(x)\qquad(x\in F^\times).
\]
\end{lemma}

\begin{proof}
Since $z=b/a$, we have $ax+b=b\left(1+\frac{x}{z}\right)$. Applying \cref{lem:Q-equivariance} gives
\[
  f(ax+b)
  =\chi(b)f\left(1+\frac{x}{z}\right)
  =\chi(b)h(x/z)=\chi(b)h_z(x).
\]
\end{proof}

Let $G=F^\times$ and let $\widehat G$ be its group of multiplicative
characters.  For a function $\varphi\colon G\to\C$, we use the Fourier
transform convention
\[
  \widehat\varphi(\theta)=\sum_{t\in G}\varphi(t)\theta(t),
  \qquad \theta\in\widehat G.
\]

Note that $\widehat \varphi : \widehat G \to \mathbb{C}$. As usual, we write
\[
  \bigl|\widehat{\varphi}(\theta)\bigr|^2
  :=
  \widehat{\varphi}(\theta)\,
  \overline{\widehat{\varphi}(\theta)}.
\]
Thus, $\bigl|\widehat{\varphi}(\theta)\bigr|^2$ is a nonnegative
real number, and it vanishes if and only if
$\widehat{\varphi}(\theta)=0$.

Since character values are roots of unity,
$\overline{\theta(t)}=\theta(t)^{-1}$.

In particular,
\begin{equation}
  \widehat h(\theta)
  =\sum_{\substack{u\in Q\\u\ne1}}\chi(u)\theta(u-1).
  \label{eq:hhat}
\end{equation}
Suppose now that $\theta$ is trivial on $\F_q^\times$. Then it factors through
$F^\times/\F_q^\times$.  The map
\[
  F^\times/\F_q^\times\longrightarrow Q,
  \qquad r\F_q^\times\longmapsto r^{q-1},
\]
is an isomorphism: the power map $r\mapsto r^{q-1}$ has kernel $\F_q^\times$ and image of order $q+1$, hence image $Q$.  We therefore
associate with $\theta$ the character $\lambda_\theta\in\widehat Q$ defined by
\begin{equation}
  \lambda_\theta(r^{q-1})=\theta(r),
  \qquad r\in F^\times.
  \label{eq:lambda}
\end{equation}
This is well defined: any two elements with the same $(q-1)$st power
differ by an element of $\F_q^\times$, on which $\theta$ is trivial.
It is a character because both the power map and $\theta$ are
multiplicative.

\begin{proposition}
\label{prop:fourier-nonzero}
Every multiplicative Fourier coefficient of $h$ is nonzero.  More
precisely, for $\theta\in\widehat{F^\times}$,
\[
  \bigl|\widehat h(\theta)\bigr|^2=
  \begin{cases}
    q,&\theta\restr_{\F_q^\times}\ne1,\\
    q^2,&\theta\restr_{\F_q^\times}=1
       \text{ and }\lambda_\theta=\chi,\\
    1,&\theta\restr_{\F_q^\times}=1
       \text{ and }\lambda_\theta\ne\chi.
  \end{cases}
\]
\end{proposition}

\par\medskip

\begin{proof}
Set $S_\theta=\widehat h(\theta)$.  From \eqref{eq:hhat}, the
multiplicativity of $\chi$, and
$\overline{\theta(t)}=\theta(t)^{-1}$, we obtain
\begin{align}
  |S_\theta|^2
  &=\sum_{\substack{u,v\in Q\\u,v\ne1}}
      \chi(u)\chi(v)\theta\!\left(\frac{u-1}{v-1}\right)\notag\\
  &=\sum_{w\in Q}\chi(w)
      \sum_{\substack{v\in Q\\v\ne1,\,wv\ne1}}
      \theta\!\left(\frac{wv-1}{v-1}\right),
  \label{eq:double-sum}
\end{align}
where $w=u/v$.

For given $w \in Q \setminus \{ 1 \}$, the number
\[ \frac{u-1}{v-1} = \frac{wv-1}{v-1} = w + \frac{w-1}{v-1} \]
takes $q-1$ distinct values as $v$ ranges over all elements of $Q \setminus \{ 1,w^{-1} \}$. As
\[ \Big( \frac{u-1}{v-1} \Big)^q = \frac{u^q-1}{v^q-1} = \frac{u^{-1}-1}{v^{-1}-1} = \frac{u-1}{v-1} \cdot \frac{v}{u},  \]
these are exactly the values of the set $\{ r \in F^\times \, : \, r^{q-1} = w^{-1} \}$. As $w^{-1} \in Q$ is a $(q-1)$st power, the latter set also has size $q-1$ as it is a coset of the kernel $\mathbb{F}_q^\ast$ of the $(q-1)$st power map.

When $w=1$, the ratio in \eqref{eq:double-sum} is $1$ for all
$q$ choices of $v\in Q\setminus\{1\}$.  Consequently,
\begin{equation}
  |S_\theta|^2
  =q+\sum_{\substack{w\in Q\\w\ne1}}\chi(w)
       \sum_{r^{q-1}=w^{-1}}\theta(r).
  \label{eq:S-square}
\end{equation}

Suppose first that $\theta\restr_{\F_q^\times}$ is nontrivial.  For a
fixed $w\ne1$, the solutions of $r^{q-1}=w^{-1}$ form a coset
$r_0\F_q^\times$.  Hence
\[
  \sum_{r^{q-1}=w^{-1}}\theta(r)
  =\theta(r_0)\sum_{c\in\F_q^\times}\theta(c)=0,
\]
and \eqref{eq:S-square} gives $|S_\theta|^2=q$.

Now suppose that $\theta\restr_{\F_q^\times}=1$.  By
\eqref{eq:lambda},
\[
  \sum_{r^{q-1}=w^{-1}}\theta(r)
  =(q-1)\lambda_\theta(w^{-1}).
\]
Therefore
\[ |S_\theta|^2=q+(q-1)\sum_{\substack{w\in Q\\w\ne1}} \chi(w)\lambda_\theta(w)^{-1}=q+(q-1)\left(\sum_{w\in Q}\chi(w)\lambda_\theta(w)^{-1}-1\right). \]

\par\smallskip
Character orthogonality on $Q$ makes the last complete sum equal to
$q+1$ when $\lambda_\theta=\chi$ and to $0$ otherwise.  The remaining
two values are therefore $q^2$ and $1$, respectively.
\end{proof}

\begin{lemma} \label{new1}
If $\mathcal{C} = Q - 1$ and $\mathcal{R} = \mathcal{C} \setminus \{ 0 \}$, then $\mathcal{R}$ is a relative $(q+1,q-1,q,1)$-difference set in $F^\times$ relative to $\mathbb{F}_q^\times$.
\end{lemma}
\begin{proof}
As $|F^\times| = q^2-1 = (q+1) \cdot |\mathbb{F}_q^\times|$ and $|\mathcal{R}| = |Q| - 1=q$, the first three parameters are obvious. It remains to show that every $r \in F^\times \setminus \mathbb{F}_q^\times$ can be written in a unique way as $\frac{v-1}{u-1}$, where $u$, $v$ are distinct elements of $Q \setminus \{ 1 \}$. We have already shown that if $r = \frac{v-1}{u-1}$, then $r^{q-1} = w^{-1}$, where $w = \frac{u}{v} \not= 1$. This then shows that $r \not\in \mathbb{F}_q$ and that $w$ is uniquely determined by $r$. Now, $r = \frac{u-1}{v-1} = w + \frac{w-1}{v-1}$ and so $v$ and hence also $u$ are uniquely determined by $r$. This implies that the $q(q-1)$ values $\frac{v-1}{u-1} \in F^\times \setminus \F_q^\times$ where $u,v \in Q \setminus 1$ with $u \not= v$ are mutually distinct and constitute the whole of $\F^\times \setminus F_q^\times$.
\end{proof} 

 \medskip We can give a geometric meaning to Lemma \ref{new1}. For $u\in Q\setminus\{1\}$, the point $u-1$ determines the direction of the secant line of $Q$ joining $1$ and $u$.  Under the identification
\[ F^\times/\F_q^\times\xrightarrow{\ \sim\ }Q, \qquad x\F_q^\times\longmapsto x^{q-1}, \]
this direction is represented by
\[ (u-1)^{q-1} =\frac{u^q-1}{u-1}=\frac{u^{-1}-1}{u-1} =-u^{-1}. \]
$\mathcal{R}$ contains exactly one point on every line through the origin except the tangent to $\mathcal{R}$ at the origin. The missing direction corresponds to $-1\in Q$. The ratio $r=\frac{u-1}{v-1}$ compares the directions of the two secants determined by $u,v\in Q \setminus \{ 1 \}$. By Lemma \ref{new1}, this ratio uniquely determines the two points $u$ and $v$ of the conic $Q$.  

\medskip Lemma \ref{new1} allows us to prove the following.

\begin{lemma} \label{new2}
Put
\[ C_h(r)=\sum_{t\in F^\times}h(rt)h(t) \]
for every $r \in F^\times$. Then 
\[  C_h=q\,\mathbf 1_{\{1\}}+\eta-\mathbf 1_{\F_q^\times}, \]
where $\eta$ is the multiplicative character of $F^\times$ defined by $\eta(r) :=\chi(r^{1-q})$, $r \in F^\times$.
\end{lemma}
\begin{proof}
For $r,t \in F^\times$, we have $h(rt) h(t) = f(1+rt) f(1+t) = f(u) f(v)$, and this is only nonzero whenever $u:=1+rt \in Q \setminus \{ 1 \}$ and $v = 1 + t \in Q \setminus \{ 1 \}$. Then $r = \frac{u-1}{v-1}$ and $r^{1-q} = \frac{u-1}{v-1} \cdot \frac{v^q-1}{u^q-1} = \frac{u-1}{v-1} \cdot \frac{1-v}{1-u} \cdot \frac{u}{v} = \frac{u}{v}$ and $\eta(r) = \chi(r^{1-q}) = \chi(\frac{u}{v}) = \chi(u) \chi(v)^{-1} = \chi(u) \chi(v) = f(u) f(v) = h(rt) h(t)$. If $r \not= 1$, then we know by Lemma \ref{new1} that here are unique possibilities for $u$, $v$ and $t$, implying that in the summation defining $C_h(r)$, there is only one nonzero term which is moreover equal to $\eta(r)$. For $r=1$, we know that there are $|Q \setminus \{ 1 \}| = q$ possibilities for $u$, $v$ and $t$. As $\eta$ is trivial on $\mathbb{F}_q$, we then see that the claim holds.       
\end{proof}

\medskip The multiplicative character that occurs in Lemma \ref{new2} plays an important role in the context of Proposition \ref{prop:fourier-nonzero}, more exactly to answer the question under which exact circumstances the characters $\lambda_\theta$ and $\chi$ of $Q$ coincide. Recall that $\lambda_\theta$ is not the restriction of $\theta$ to $Q$, but is transported from the quotient $F^\times/\F_q^\times$ through \eqref{eq:lambda}.

\begin{lemma} \label{new3}
There is a unique character $\theta$ of $F^\times$ for which $\lambda_\theta = \chi$ and $\theta\restr_{\F_q^\times}=1$, namely $\theta = \eta$.
\end{lemma}
\begin{proof}
As  $\eta$ is trivial on $\F_q^\times$, \eqref{eq:lambda} implies that there is (at most) one such character. Put $F^\times = \langle g \rangle$. Then $\F_q^\times=\langle g^{q+1}\rangle$ and $Q=\langle g^{q-1}\rangle$. The transported character of $\chi$ satisfies $\lambda_\eta(g^{q-1})=\eta(g)=\chi(g^{1-q})=\chi(g^{q-1})$ and hence $\lambda_\eta=\chi$.
\end{proof}

\bigskip Lemma \ref{new2} also allows to give an alternative proof of Proposition \ref{prop:fourier-nonzero}. Indeed, by invoking that lemma, we have that
\[ \bigl|\widehat h(\theta)\bigr|^2 =  \widehat h(\theta) \, \overline{\widehat h(\theta)} =   \sum_{u,t\in F^\times} h(u)h(t)\theta(u)\theta(t)^{-1} =   \sum_{r,t\in F^\times}h(rt)h(t)\theta(r) =  \sum_{r\in F^\times}C_h(r)\theta(r)  \]
is equal to $q+\sum_{r\in F^\times}\eta(r)\theta(r)-\sum_{r\in\F_q^\times}\theta(r)$. By character orthogonality,
\[
  \sum_{r\in F^\times}\eta(r)\theta(r)
  =
  \begin{cases}
    q^2-1,&\theta=\eta,\\
    0,&\theta\ne\eta,
  \end{cases}
\]
whereas
\[
  \sum_{r\in\F_q^\times}\theta(r)
  =
  \begin{cases}
    q-1,&\theta\restr_{\F_q^\times}=1,\\
    0,&\theta\restr_{\F_q^\times}\ne1.
  \end{cases}
\]
Substitution gives the respective values $q$, $q^2$, and $1$ in \cref{prop:fourier-nonzero} (recall also Lemma \ref{new3}).

\medskip The nonvanishing in \cref{prop:fourier-nonzero} has an immediate linear-algebraic consequence.

\begin{theorem}
\label{thm:all-translates}
The family of all multiplicative translates
\[
  \{h_z:z\in F^\times\},\qquad h_z(x)=h(x/z),
\]
is a basis of $\C^{F^\times}$.
\end{theorem}

\begin{proof}
Suppose that
\[
  \sum_{z\in F^\times}c_z h(x/z)=0
  \qquad(x\in F^\times).
\]
Multiply by $\theta(x)$ and sum over $x\in F^\times$.  The substitution
$t=x/z$ gives
\[ 0 = \sum_{x \in F^\times} \sum_{z \in F^\times} c_z \theta(x) h(\frac{x}{z}) = \sum_{t \in F^\times} \sum_{z \in F^\times} c_z \theta(tz) h(t) \]
\[  = \sum_{t \in F^\times} \sum_{z \in F^\times} c_z \theta(z) \theta(t) h(t) = \widehat{h}(\theta) \sum_{z \in F^\times} c_z \theta(z).  \]
By \cref{prop:fourier-nonzero}, every factor $\widehat h(\theta)$ is nonzero.  Therefore
\begin{equation}
  \sum_{z\in F^\times}c_z\theta(z)=0
  \qquad\text{for every }\theta\in\widehat{F^\times}.
  \label{eq:character-relations}
\end{equation}

For completeness, we spell out the nonsingularity argument.  Choose a
generator $g$ of $F^\times$, put $M=q^2-1$, and let
$\omega=\exp(2\pi i/M)$.  The characters of $F^\times$ are
\[
  \theta_j(g^m)=\omega^{mj},
  \qquad 0\le j,m\le M-1.
\]
Writing \eqref{eq:character-relations} for these $M$ characters gives
\[
  \sum_{m=0}^{M-1}c_{g^m}\omega^{mj}=0,
  \qquad 0\le j\le M-1.
\]
The coefficient matrix is the Vandermonde matrix
\[
  V=(\omega^{mj})_{0\le j,m\le M-1},
\]
whose determinant is
\[
  \det V=\prod_{0\le r<s\le M-1}(\omega^s-\omega^r)\ne0.
\]
Thus $c_{g^m}=0$ for every $m$, and the translates are linearly
independent.  There are $M=|F^\times|$ of them in the
$M$-dimensional space $\C^{F^\times}$, so they form a basis.
\end{proof}

{
The preceding theorem is the main Fourier statement.  We now record
its consequences for the two nonprincipal eigenspaces.  Let
$D=(F^\times)^2$, let $\Gamma=P(q^2)$, and write
\[
  r=\frac{q-1}{2},\qquad s=-\frac{q+1}{2}
\]
for the two nonprincipal eigenvalues.  Their eigenspaces $E_r$ and
$E_s$ both have dimension $(q^2-1)/2$.

We first recall the effect of a nonsquare dilation.  The adjacency
operator of $\Gamma$ is
\[
  (\mathsf A g)(x)=\sum_{d\in D}g(x+d).
\]

\begin{proposition} \label{prop:switch}
Let $g$ be a nonprincipal eigenfunction of $P(q^2)$ with eigenvalue $\lambda$. Then:
\begin{itemize}
\item if $\mu \in F^\times$ is a square, then $g_\mu(x) = g(\mu x)$ is an eigenfunction with eigenvalue $\lambda$.
\item if $\mu \in F^\times$ is a nonsquare, then $g_\mu(x) = g(\mu x)$ is an eigenfunction with eigenvalue $-1-\lambda$.
\end{itemize}
\end{proposition}

\begin{proof}
Every nonprincipal eigenfunction is orthogonal to the all-one vector, so $\sum_{y\in F}g(y)=0$. Consequently, for every $x\in F$,
\begin{equation}
  \sum_{e\in F^\times}g(\mu x+e)=-g(\mu x).
  \label{eq:sum-nonzero-shifts}
\end{equation}
By the definition of the adjacency operator and of $g_\mu$, we have
\[
  (\mathsf A g_\mu)(x)=\sum_{d\in D}g_\mu(x+d)
  =\sum_{d\in D}g\bigl(\mu(x+d)\bigr)=\sum_{d\in D}g(\mu x+\mu d).
\]
If $\mu$ is a square, then
\[   (\mathsf A g_\mu)(x) = \sum_{d \in D} g(\mu x + \mu d) = \sum_{e \in D} g(\mu x + e) = \lambda g(\mu x) = \lambda g_\mu(x).  \]
If $\mu$ is a nonsquare, then $\mu D=F^\times\setminus D$ and we compute
\begin{align*}
  (\mathsf A g_\mu)(x)
  &=\sum_{d\in D}g(\mu x+\mu d)
    =\sum_{e\in F^\times\setminus D}g(\mu x+e)\\
  &=\sum_{e\in F^\times}g(\mu x+e)
    -\sum_{e\in D}g(\mu x+e)\\
  &=-g(\mu x)-\lambda g(\mu x)
    =(-1-\lambda)g_\mu(x).
\end{align*}
Here, we have used \eqref{eq:sum-nonzero-shifts} and the eigenfunction equation for $g$ at $\mu x$.
\end{proof}

\begin{corollary}
\label{cor:z-eigenspace-bases}
For $z\in F^\times$, extend $h_z$ from $F^\times$ to $F$ by
\[
  \widetilde h_z(x)=f\left(1+\frac{x}{z}\right).
\]
Then every $\widetilde h_z$ is an optimal eigenfunction with support
of size $q+1$.  More precisely,  $\{\widetilde h_z:z\in D\}$ is a basis of the same nonprincipal eigenspace as $f$, whereas $\{\widetilde h_z:z\in F^\times\setminus D\}$
is a basis of the other nonprincipal eigenspace.  
The full family $\{\widetilde h_z:z\in F^\times\}$ is therefore a
basis of $E_r\oplus E_s=\mathbf1^\perp$.
\end{corollary}

\begin{proof}
The restriction of $\widetilde h_z$ to $F^\times$ is $h_z$.
Therefore \cref{thm:all-translates} implies that the $q^2-1$
functions $\widetilde h_z$ are linearly independent.

The function $\widetilde h_z$ is obtained from $f$ by the invertible
affine transformation
\[
  x\longmapsto z^{-1}x+1.
\]
In particular,
\[
  \supp(\widetilde h_z)=z(Q-1),
\]
so its support has size $q+1$.  If $z\in D$, then $z^{-1}$ is a
square and the displayed affine transformation is an automorphism of
the Paley graph.  Hence $\widetilde h_z$ is an optimal eigenfunction
in the same eigenspace as $f$.

If $z\notin D$, translation by $1$ preserves the eigenspace of $f$,
while dilation by the nonsquare $z^{-1}$ sends its eigenvalue
$\lambda$ to $-1-\lambda$ by \cref{prop:switch}.  Thus
$\widetilde h_z$ belongs to the other nonprincipal eigenspace.  This
dilation is an invertible and support-preserving correspondence between
the two eigenspaces.  If the second eigenspace contained an
eigenfunction with support smaller than $q+1$, applying the inverse
nonsquare dilation would contradict the optimality of $f$.  Hence the
nonsquare-indexed functions are optimal as well.

Each of $D$ and $F^\times\setminus D$ has $(q^2-1)/2$ elements.
The corresponding functions are independent and lie in an eigenspace
of that dimension, so each half is a basis.  The assignment of the
eigenvalues follows from \cref{prop:optimal}.  Taking both halves gives
the asserted basis of $E_r\oplus E_s=\mathbf1^\perp$.
\end{proof}

\medskip \noindent The following is a consequence of Proposition \ref{prop:switch} and Corollary \ref{cor:z-eigenspace-bases}.

\begin{corollary}
\label{cor:affine-form}
Fix $b\in Q$.  Then
\[
  \bigl\{x\longmapsto f(ax+b):a\in D\bigr\}
\]
is a basis of the same nonprincipal eigenspace as $f$, whereas
\[
  \bigl\{x\longmapsto f(ax+b):a\in F^\times\setminus D\bigr\}
\]
is a basis of the other nonprincipal eigenspace.  Every function in
the two families is optimal.  Their union $\bigl\{x\longmapsto f(ax+b):a\in F^\times\bigr\}$ is a basis of $E_r\oplus E_s=\mathbf1^\perp$.
\end{corollary}

\begin{remark}
The distinction between $h_z$ and $\widetilde h_z$ is important.
Optimality refers to eigenfunctions on the full vertex set $F$, so it
applies to $\widetilde h_z$, not directly to its restriction $h_z$.
Since $\widetilde h_z(0)=f(1)=1$, one has
\[
  |\supp(\widetilde h_z)|=q+1,
  \qquad
  |\supp(h_z)|=q.
\]
\end{remark}

We can now complete the proofs of the main results stated in the introduction.

\begin{proof}[Proof of \cref{thm:main}]
The assertion on the restrictions $h_z$ is
\cref{thm:all-translates}.  All statements concerning their affine
extensions $\widetilde h_z$, including optimality, the square and
nonsquare eigenspace bases, and the basis of $\mathbf1^\perp$, are
proved in \cref{cor:z-eigenspace-bases}.  Finally, the functions are
real-valued, so their complex linear independence also gives the
corresponding assertions over $\mathbb R$.
\end{proof}

\begin{proof}[Proof of \cref{co:main}]
For a function $g: F^\times \to \mathbb{R}$ to belong to a certain eigenspace, it must be orthogonal to the other eigenspaces. For the eigenspace corresponding to the valency of $P(q^2)$, this amounts to condition (P1). For the other eigenspace, this amounts to the other mentioned condition (P$i$). This is a consequence of Corollary \ref{cor:z-eigenspace-bases}, after expressing what orthogonality with the eigenfunction $\widetilde h_z$, $z \in F^\times$, means. Indeed, the map $g$ is orthogonal with such an eigenfunction if
\[ 0 = \sum_{x \in F} g(x) \widetilde{h}_z(x) = \sum_{x \in F} g(x) \cdot f(1 + \frac{x}{z}) \]
\[ = \sum_{x \in z (Q-1)} g(x) f(1 + \frac{x}{z}) = \sum_{x \in z(\Qzero-1)} g(x) - \sum_{x \in z(\Qone-1)} g(x). \]
\end{proof}

\section{Examples} \label{sec:examples}

We illustrate the construction in the smallest field from each congruence class modulo 4, emphasizing how square and nonsquare values of $z$
produce the two eigenspace bases.

\begin{example}[The case $q\equiv1\pmod4$] \label{ex:q5}
Let $q=5$. The graph $P(25)$ has $25$ vertices, valency $12$, and nonprincipal eigenvalues $2$ and $-3$. The norm-one group $Q$ has order $6$, while $\Qzero$ and $\Qone$ have order $3$. Since $5 \equiv 1 \pmod4$, the graph induced on $Q$ is $K_{3,3}$ with bipartition $Q=\Qzero\sqcup\Qone$. The function $f$, equal to $1$ on one part and $-1$ on the other, is the optimal eigenfunction for the eigenvalue $-3$ described in \cref{prop:optimal}.

For a concrete realization, take $\F_{25}=\F_5(\alpha)$, where $\alpha^2=2$. Then $N(X+Y\alpha)=X^2-2Y^2$, and a direct calculation gives
\begin{align*}
  Q&=\{1,-1,2+2\alpha,2-2\alpha,
                 3+2\alpha,3-2\alpha\},\\
  \Qzero&=\{1,2+2\alpha,2-2\alpha\},\\
  \Qone&=\{-1,3+2\alpha,3-2\alpha\}.
\end{align*}
Let $D=(\F_{25}^\times)^2$.  The square-indexed functions
\[
  \left\{\widetilde h_z(x)=f\left(1+\frac{x}{z}\right):z\in D\right\}
\]
form a basis of the $(-3)$-eigenspace, while the nonsquare-indexed
functions
\[
  \left\{\widetilde h_z(x)=f\left(1+\frac{x}{z}\right):
          z\in\F_{25}^\times\setminus D\right\}
\]
form a basis of the $2$-eigenspace.  Each family contains $12$ optimal eigenfunctions.

The twelve square parameters are precisely $z \in \{1,-1,2,-2,1+\alpha,1-\alpha,-1+\alpha,-1-\alpha,2+2\alpha,2-2\alpha,-2+2\alpha,-2-2\alpha \}$, while the non-square parameters are given by $z \in \{ \alpha,-\alpha,2\alpha,-2\alpha,1+2\alpha,1-2\alpha,-1+2\alpha,-1-2\alpha,2+\alpha,2-\alpha,-2+\alpha,-2-\alpha \}$.
\end{example}

\begin{example}[The case $q\equiv3\pmod4$]
\label{ex:q3}
Let $q=3$. The graph $P(9)$ has $9$ vertices, valency $4$, and nonprincipal eigenvalues $1$ and $-2$. Here $|Q|=4$ and $|\Qzero|=|\Qone|=2$.  Since
$3\equiv3\pmod4$, the graph induced on $Q$ is $K_2\mathbin{\dot\cup}K_2$, with components induced by $\Qzero$ and $\Qone$.  The function $f$ is an optimal eigenfunction for the eigenvalue $1$.

To make the construction explicit, write $\F_9=\F_3(\alpha)$, where $\alpha^2=-1$. Here $N(X+Y\alpha)=X^2+Y^2$, and
 \[
  Q=\{1,-1,\alpha,-\alpha\},\qquad
  \Qzero=\{1,-1\},\qquad
  \Qone=\{\alpha,-\alpha\}.
\]
Let $D=(\F_9^\times)^2$.  The four square-indexed functions
\[
  \left\{\widetilde h_z(x)=f\left(1+\frac{x}{z}\right):z\in D\right\}
\]
form a basis of the $1$-eigenspace, whereas the four nonsquare-indexed functions
\[
  \left\{\widetilde h_z(x)=f\left(1+\frac{x}{z}\right):
          z\in\F_9^\times\setminus D\right\}
\]
form a basis of the $(-2)$-eigenspace.

In this example $Q=D$, since both are the unique subgroup of order $4$ in $\F_9^\times$.  Thus the four square parameters are precisely $z\in\{1,-1,\alpha,-\alpha\}$, while the non-square parameters are given by $z \in \{ 1+\alpha,1-\alpha,-1+\alpha,-1-\alpha \}$.
\end{example}

\section*{Acknowledgments}
The first four authors thank the Dipartimento di Matematica e Fisica of the Universit\`a degli Studi della Campania ``Luigi Vanvitelli'' for the hospitality enjoyed during a visit in June 2026 during which this research was initiated. The last author gratefully acknowledges the hospitality of the School of Mathematical Sciences at Hebei Normal University, where part of this research was carried out during July and August 2026. 
The second author is supported by the Natural Science Foundation of Hebei Province (A2023205045) and the 111 Center (Grant No.D26018). The third author is supported by Science Foundation of Hebei Normal University (L2026B48).
The fourth author is supported by China Scholarship Council (202608130044).
Also, the last author was partially supported by the Italian National Group for Algebraic and Geometric Structures and their Applications (GNSAGA - INdAM).

\end{document}